\documentclass[12pt]{article}
\usepackage{Iopconf}
\usepackage{times,colordvi,amsmath,epsfig,float,color,subfigure}
\usepackage{amssymb}
\usepackage{amsthm}
\usepackage{amsopn}
\usepackage{amscd}
\usepackage[latin1]{inputenc}
\input xy
\xyoption{all}

\newtheorem{thm}{Theorem}

\newtheorem{defn}[thm]{Definition}

\newcommand{\be}{\begin{equation}} \newcommand{\ee}{\end{equation}}
\newcommand{\ba}{\begin{eqnarray}} \newcommand{\ea}{\end{eqnarray}}

\def\e{\epsilon}
\def\ve{\varepsilon}

\def\s{\sigma}
\def\D{\Delta}

\def\tR{{\tilde{\cal R}}}

\def\C{\cal C}
\def\A{\cal A}

\def\F{\cal F}
\def\tp{\otimes}

\def\I{{\rm id}}

\def\rar{\rightarrow}

\begin{document}

\title{Ribbon structure in symmetric pre-monoidal categories}

\author{Liam Wagner\footnote{ldw@maths.uq.edu.au}, Jon 
Links\footnote{jrl@maths.uq.edu.au} and 
Phil Isaac\footnote{psi@maths.uq.edu.au}}

\affil{Centre for Mathematical Physics, The University of
     Queensland, Brisbane, 4072, Australia}

\beginabstract
Let $U(g)$ denote the universal enveloping algebra of a Lie algebra $g$. 
We show the existence of a ribbon algebra structure in a 
particular deformation of $U(g)$ which leads to a symmetric 
pre-monoidal category of $U(g)$-modules.  
\endabstract
\section{Introduction} 

Recently there has been some work \cite{y,joyce1,joyce2,joyce3} on
investigating the properties of categories in which the pentagon axiom
which is central to the notion of monoidal categories \cite{mac} is
forsaken. In particular it was shown in \cite{isaaci} that a certain
deformation of the Hopf algebra structure of the universal enveloping algebras
of
Lie algebras naturally gives rise to a construction for symmetric
pre-monoidal categories as defined in \cite{joyce1,joyce2,joyce3}. 
It is then reasonable to ask
to what extent can such categories also be endowed with a ribbon
structure, in analogy with the known examples of ribbon categories arising
from representations of ribbon quasi-Hopf algebras \cite{ac}. 
Here we note that in a
particular specialisation of the deformation given in \cite{isaaci} one can
indeed preserve the ribbon algebra structure which provides an initial
insight into making the notions of ribbon and pre-monoidal categories
compatible.
      
\section{Symmetric pre-monoidal categories}

We shall begin by defining a pre-monoidal category, taken from
\cite{isaaci}.  \begin{defn} A pre-monoidal category is a triple $({\C} ,
\tp, \A)$ where $\C$ is the category of objects, $\tp$ is a bifunctor $\tp
: \C \tp \C \rar \C$ and $\A$ is the natural associator isomorphism such
that ${\A}: \tp(\I \times \tp) \rar \tp (\tp \times \I)$. \end{defn}

\noindent 
In practice we we write $a_{U,V,W}$ as the action of the associator such that 
$$ a_{U,V,W}: U\tp (V \tp W) \rightarrow (U \tp V) \tp W.$$ 
We note that there are no conditions imposed on $\A$,
but it is important that we should define the natural isomorphism 
$Q:\tp (\tp \times \tp) \rar \tp (\tp \times \tp)$ via the following diagram:

$$
\xymatrix{
(U\tp V)\tp(W\tp Z)  \ar[dd]_{a_{(U\tp V),W,Z}}                           &
                                                                        & 
(U\tp V)\tp(W\tp Z) \ar[ll]_{q_{U,V,W,Z}}                               \\ &&\\
((U\tp V)\tp W)\tp Z                                                    &
                                                                        &
U\tp (V\tp (W\tp Z)) \ar[uu]_{a_{U,V,(W\tp Z)}} \ar[dd]^{\I \tp a_{V,W,Z}} \\
&&\\
(U\tp(V\tp W))\tp Z \ar[uu]^{a_{U,V,W}\tp \I}                           &
                                                                        &
U\tp ((V\tp W)\tp Z) \ar[ll]^{a_{U,(V\tp W),Z}}    
}
$$

\noindent This box diagram expresses $Q$ in terms of the 
associator isomorphisms which are
defined as:
\begin{equation}
Q((U\tp V)\tp (W \tp Z))= q_{U,V,W,Z}((U \tp V) \tp (W \tp Z)).  
\end{equation}
\noindent These may be expressed  as
\begin{equation}
q_{U,V,W,Z} =a^{-1}_{(U\tp V),W,Z}(a_{U,V,W} \tp \I)a_{U,(V \tp W),Z}(\I \tp
a_{V,W,Z})a^{-1}_{U,V,(W\tp Z)}.
\end{equation}

This condition is a generalisation of the pentagon condition 
used in describing monoidal categories, where $Q=\I \tp \I \tp \I \tp \I$. We can examine the significance of $Q$ and its
use in distinguishing the coupling of the objects of the category through 
employing distinct brackets
$[]$, $\{ \}$. This notation shows that 
\begin{equation}
Q([U\tp V] \tp \{W \tp Z\}) = (\{ U \tp V \} \tp [W \tp Z]).  
\end{equation}

We describe pre-monoidal categories as being unital if they have 
an identity object ${\mathbf 1}\in \C$ and natural
isomorphisms $\rho_U: U\otimes {\mathbf 1}\rightarrow U$ and $\lambda_U:
{\mathbf 1}\otimes U\rightarrow U$. The next important class of unital pre-monoidal categories is when the tensor product is
commutative up to isomorphism. This leads to the concept of a braided unital pre-monoidal
category. 

\begin{defn}\label{bmc}
A unital pre-monoidal category $\C$ is said to be braided if it is equipped with
a natural commutativity isomorphism $\s_{U,V}:U\tp V \rar V\tp
U$ for all objects $U,V\in \C$ such that the following diagrams
commute \cite{joyce3}:
\begin{itemize}
\item[(i)] 
$$
\xymatrix{
                                                                & 
U\tp (V\tp W) \ar[rr]^{\s_{U,(V\tp W)}}                           &&
(V\tp W)\tp U \ar[dr]^{a^{-1}_{V,W,U}}                         &               
                                                                \\
(U\tp V)\tp W \ar[ur]^{a^{-1}_{U,V,W}} \ar[dr]_{\s_{U,V}\tp\I} &&&&
V\tp (W\tp U)                                                   \\
                                                                &
(V\tp U)\tp W \ar[rr]_{a^{-1}_{V,U,W}}                         &&
V\tp (U\tp W) \ar[ur]_{\I\tp\s_{U,W}}                           &
}
$$ 
\item[(ii)]
$$
\xymatrix{
                                                           & 
(U\tp V)\tp W \ar[rr]^{\s_{(U\tp V),W}}                      &&
W\tp (U\tp V) \ar[dr]^{a_{W,U,V}}                         &               
                                                           \\
U\tp (V\tp W) \ar[ur]^{a_{U,V,W}} \ar[dr]_{\I\tp\s_{V,W}} &&&&
(W\tp U)\tp V                                              \\
                                                           &
U\tp (W\tp V) \ar[rr]_{a_{U,W,V}}                         &&
(U\tp W)\tp V \ar[ur]_{\s_{U,W}\tp\I}                      &
}
$$
\item[(iii)]
$$
\xymatrix{
(U\tp V)\tp (W\tp Z) \ar[d]_{\s_{(U\tp V),(W\tp Z)}} \ar[rr]^{q_{U,V,W,Z}} &&
(U\tp V)\tp (W\tp Z) \ar[d]^{\s_{(U\tp V),(W\tp Z)}} \\
(W\tp Z)\tp (U\tp V) && (W\tp Z)\tp (U\tp V) \ar[ll]_{q_{W,Z,U,V}}
}
$$
\end{itemize}
\end{defn}

Hereafter we will only be concerned with the case of {\it symmetric} categories in which the commutativity isomorphism satisfies the additional condition $\sigma_{U,V}\circ \sigma_{V,U} =
\I_{V \tp U}$ for all objects $U,\,V \in \C$, and $\circ$ denotes composition of morphisms. In this case the diagrams ($i$) and ($ii$) are equivalent.  

\section{Twining}
It is known that the finite-dimensional modules of quasi-triangular 
quasi-bialgebras give rise to braided
monoidal categories, and that these categories are invariant under twisting 
\cite{drinfeld,cp,turaev}. Here we need to use a different method to build 
pre-monoidal categories. We follow the work of \cite{isaaci} where a twining operation was introduced to achieve this.
\\\\
Define $A$ as a quasi-triangular quasi-bialgebra by the octuple 
$(A,\D,\ve , {\Phi}, {\mathcal{R}},S,\alpha,\beta)$ where 
$\D: A \rightarrow A \tp A$ is the co-product, $\Phi$ is an invertible element
called the co-associator which is defined as $\Phi= I \tp I \tp I$, 
$\ve : A \rightarrow \mathbb{C}$ is the co-unit, $\mathcal{R} \in A \tp A$ is the universal $R$-matrix satisfying 
$$\mathcal{R} \D(x) = \D^T(x) \mathcal{R} ~~~~~\forall\,x\in A $$ 
where $\D^T$ is the opposite co-product. Further $S$ is the antipode 
and $\alpha,\,\beta$ are canonical elements satisfying certain properties 
(see \cite{drinfeld,cp,turaev} for details).  
Let $K$ be an element of the centre of $A$ and let
non-zero $\gamma \in \mathbb{C}$ be fixed but arbitrary. 
Following \cite{isaaci} we define 
\begin{equation} 
\begin{split}
&\tilde{\mathcal{R}} = \gamma ^{K \tp K} \mathcal{R}= \mathcal{R}\cdot
\gamma^{K\tp K}\\
&\tilde{\Phi}=\Phi \cdot \gamma^{\kappa} 
\end{split}
\label{twining}
\end{equation}

\noindent such that $\kappa =K\tp (I \tp K + K\tp I - \D(K))$. 
In general  we write 
\begin{equation}
\begin{split}
\tilde{\Phi} &= \sum_{j} X_{j} \tp Y_{j} \tp Z_{j}\\
{\tilde{{\Phi}}^{-1}} &= \sum_{j} \bar{X_{j}} \tp \bar{Y_{j}} \tp \bar{Z_{j}}\\
\tilde{\mathcal{R}} &= \sum_{j} a_{j} \tp b_{j}.
\end{split}
\label{randphi}
\end{equation}
Then the following  hold:

\begin{equation}
\begin{split}
&(\I \tp \D) \D (a) =\tilde{\Phi}^{-1}(\D \tp \I)\D(a) \tilde{\Phi} ~~\; \; \forall
\,a \in A, \\
& \tilde{\mathcal{R}} \D (a) = \D^{T}(a) \tilde{\mathcal{R}},\\
&(\D \tp \I) \tilde{\mathcal{R}} = \tilde{\Phi}_{312}\tilde{\mathcal{R}}_{13} \tilde{\Phi}_{132}
\tilde{\mathcal{R}}_{23} \tilde{\Phi}^{-1}_{123},\\
&(\I \tp \D) \tilde{\mathcal{R}} = \tilde{\Phi}_{213}^{-1}\tilde{\mathcal{R}}_{13} \tilde{\Phi}_{213}^{-1}
\tilde{\mathcal{R}}_{12} \tilde{\Phi}_{123}(\gamma^{2\kappa})^{-1}_{123}
\end{split}
\end{equation}

\noindent and we go on to define
\begin{equation}
\xi = (\D \tp \I \tp \I) \tilde{\Phi}^{-1} (\tilde{\Phi} \tp I)\cdot (\I \tp \D \tp \I) \tilde{\Phi}
\cdot ( I \tp \tilde{\Phi}) \cdot (\I \tp \I \tp \D) \tilde{\Phi}^{-1}.
\end{equation}

\noindent These relations show that the category ${\mathbf{mod}}_K(A)$ of
$A$-modules with 
\begin{equation}
a_{U,V,W} = (\pi _{U} \tp \pi _{V} \tp \pi _{W} )\tilde{\Phi},
\end{equation}

\noindent is a pre-monoidal category as $\tilde{\Phi}$ fails the pentagon condition. In particular the representation 
\begin{equation}
q_{U,V,W,Z}=(\pi_{U} \tp \pi_{V} \tp \pi_{Z})\xi
\end{equation}

\noindent does not act as the identity.

We should note that $\tR$ cannot be used in general as an ingredient in the construction of a 
braided, pre-monoidal category of $A$-modules, as the $\gamma^{2\kappa}$ defined above in 
(\ref{twining}) leads to a violation of the hexagon condition $(i)$ of 
Definition \ref{bmc}. However, it was observed in \cite{isaaci} that in the case of the universal enveloping algebra $U(g)$ of a Lie algebra $g$, it
is possible to choose  $\gamma=-1$ and $K\in\,A$ taking integer eigenvalues on all irreducible finite-dimensional $U(g)$-modules such that  ${\mathbf{mod}}_K(A)$  does possess the structure of a symmetric pre-monoidal category. 
For a given central element $K\in\,U(g)$ such that
$\ve(K)=0$ we call  $(U(g),\D,\ve,\tilde{\Phi},\tilde{\mathcal{R}})$ 
with $\tilde{\Phi}$ and $\tilde{\mathcal{R}}$ given  by (\ref{twining}) a {\it twining} of $U(g)$.

\section{Ribbon algebra structure} 

In analogy with the definition of ribbon quasi-Hopf algebras \cite{ac}, 
we can also investigate the ribbon structure for twined algebras. 
First recall that each quasi-triangular quasi-Hopf algebra $A$ possesses 
a distinguished invertible element $u$ satisfying 
\begin{eqnarray}
S^2(a)&=&uau^{-1} ~~~~~~\forall\, a\in\,A, \label{ssquared} \\
S(\alpha) u&=& \sum_j S(b_j) \alpha a_j. \label{a} \end{eqnarray} 
We then have from \cite{ac}:

\begin{defn} 
Let $A$ be a quasi-triangular quasi-Hopf algebra. We say that $A$  is a ribbon quasi-Hopf algebra if there exists a central element $v\in\,A$ such that
\begin{itemize} 
\item[i.] $v^2= u S(u)$ 
\item[ii.] $S(v)=v $
\item[iii.] $\ve (v)=1 $
\item[iv.] $\D(u v^{-1})= {\F}^{-1}(S \tp S) {\F}_{21} (uv^{-1} \tp uv^{-1})$  
\end{itemize}
where $\F$ is the Drinfeld twist \cite{drinfeld} defined by the condition
 $$\D(a)={\F}^{-1}((S\tp S)\D^T (S^{-1}(a))){\F} ~~~~~\forall\,a \in A. $$
\end{defn}

The universal enveloping algebra
$U(g)$ of a Lie algebra $g$ acquires the structure of a quasi-bialgebra with mappings:

\begin{equation}
\begin{split} 
&\ve(I)=1, \;\ve(x)=0,\; \forall x\in g \\
&S(I)=I, \; S(x)=-x, \; \forall x\in g \\
&\D(I)=I\tp I, \; \D(x)=I\otimes x+x\tp I, \; \forall x\in g  
\end{split}
\end{equation}
that are extended to all $U(g)$ such that $\e$ and $\D$ are algebra
homomorphisms and $S$ is an anti-automorphism. 
It is easily checked that $\D$ is co-associative; i.e. 
$$(\I\tp \D)\D(x)=(\D \tp \I)\D(x) ~~~~~ \forall x\in U(g) $$ 
and co-commutative 
$$\D(x) = \D^T(x) ~~~~~~\forall x\in\, U(g) $$ 
and that $S^2=\I$. This means that we can equip $U(g)$ with the structure of a quasi-triangular 
quasi-Hopf algebra by taking $\Phi=I\tp I\otimes I$ for the co-associator of $U(g)$, $\mathcal{R}=I \tp I$ as the universal $R$-matrix and $\alpha=\beta=I$. Note that in this instance the Drinfeld twist of $U(g)$ is trivial, and $U(g)$ trivially satisfies the conditions of a ribbon Hopf algebra with the choice 
$u=v=I$. 

Under twining by a central element $K$ satisfying $S(K)=-K$ we then have from (\ref{randphi})  
\begin{eqnarray} 
\tilde{\mathcal{R}}&=& (-1)^{K\otimes K} \nonumber \\
\tilde{\Phi}&=&(-1)^\kappa \nonumber 
\end{eqnarray}
and moreover it can be shown that $\alpha=\beta=I$ (cf. \cite{isaaci}). It is easily verified that the choice 
$$v=u=(-1)^{-K^2} $$ 
satisfies equations (\ref{ssquared}) and (\ref{a}) as well as the conditions 
of Definition 3. 
Thus in this instance we can conclude that the twined algebra 
$U(g)$ can still be considered as a ribbon algebra.     

\section{Conclusion}

We have shown that a particular case of the twining deformation described
in \cite{isaaci} is compatible with the notion of a ribbon structure which
can be endowed on quasi-triangular quasi-Hopf algebras. In future work we
will explore the compatibility of ribbon categories and symmetric
pre-monoidal categories on a general level and in particular investigate
the consequences of this for defining traces and inner products on generic
symmetric pre-monoidal categories (cf. \cite{kirillov}).

\end{document}